\documentclass{amsart}[11pt,amssymb]
\DeclareMathSymbol{\twoheadrightarrow} {\mathrel}{AMSa}{"10}

\def\Q{{\mathbf Q}}

\def\Z{{\mathbf Z}}
\def\C{{\mathbf C}}

\def\H{{\mathrm H}}

\def\End{\mathrm{End}}

\def\Aut{\mathrm{Aut}}
\def\Hom{\mathrm{Hom}}

\def\cl{\mathrm{cl}}

        \def\K_a{\bar{K}}

\def\OC{{\mathcal O}}
\def\P{{\mathbf P}}

\def\m{{\mathfrak m}}

\newtheorem{thm}{Theorem}[section]
\newtheorem{lem}[thm]{Lemma}
\newtheorem{cor}[thm]{Corollary}

\theoremstyle{definition}

\title[Cubic surfaces and intermediate jacobians]{Cubic surfaces and cubic threefolds, jacobians and intermediate jacobians }
\author{Yuri G. Zarhin}
\address{Department of Mathematics, Pennsylvania
State University, University Park, PA 16802, USA}
\email{zarhin\char`\@math.psu.edu}

\address{Steklov Mathematical Institute, Russian Academy of Sciences, Moscow, Russia}

\dedicatory{To my teacher Yuri Ivanovich Manin with admiration and gratitude}

\begin{document}
\maketitle

In this paper we study principally polarized abelian varieties that admit an automorphism of
order 3. It turns out that certain natural conditions on the multiplicities of its action on
the differentials of the first kind do guarantee that those polarized varieties are not
jacobians of curves. As an application, we get another proof of the (already known)
fact that intermediate jacobians of certain cubic threefolds are not jacobians of curves.

\section{Principally polarized abelian varieties that admit an automorphism of order $3$}
Let $\zeta_3=\frac{-1+\sqrt{-3}}{2}$ be a primitive (complex) cubic root of unity.
It generates  the multiplicative order $3$ cyclic group $\mu_3$ of cubic roots of unity.

Let $g>1$ be an  integer and  $(X,\lambda)$ a principally polarized $g$-dimensional
abelian variety over the field $\C$ of complex numbers, $\delta$ an automorphism of
$(X,\lambda)$ that satisfies the cyclotomic equation
$\delta^2+\delta+1=0$ in $\End(X)$. In other words, $\delta$ is a periodic automorphism
of order $3$, whose set of fixed points is finite. This gives rise to the embeddings
$$\Z[\zeta_3]\hookrightarrow \End(X), 1\mapsto 1_X,\ \zeta_3\mapsto \delta,$$
$$\Q(\zeta_3)\hookrightarrow \End^{0}(X), 1\mapsto 1_X,\ \zeta_3\mapsto \delta.$$
By functoriality, $\Q(\zeta_3)$ acts on the $g$-dimensional complex
vector space $\Omega^1(X)$ of  differentials of the first kind on
$X$. This provides $\Omega^1(X)$ with a structure of
 $\Q(\zeta_3)\otimes_{\Q}\C$-module. Clearly,
 $$\Q(\zeta_3)\otimes_{\Q}\C=\C\oplus \C$$
 where the summands correspond to the embeddings $\Q(\zeta_3) \to \C$ that send
 $\zeta_3$ to $\zeta_3$ and $\zeta_3^{-1}$ respectively. So, $\Q(\zeta_3)$ acts on
 $\Omega^1(X)$ with multiplicities $a$ and $b$ that correspond to the two embeddings of
  $\Q(\zeta_3)$ into $\C$. Clearly,
 $a$ and $b$ are non-negative integers with $a+b=g$.

 \begin{thm}
 \label{nonjacob}
If $g+2<3\mid a-b\mid$ then $(X,\lambda)$ is not the jacobian of a smooth projective irreducible
genus $g$ curve with canonical principal
 polarization.
 \end{thm}

 \begin{proof}
Suppose that $(X,\lambda)\cong (J(C),\Theta)$ where $C$ is an
irreducible smooth projective genus $g$ curve, $J(C)$ its jacobian
with canonical principal polarization $\Theta$. It follows from the
Torelli theorem in Weil's form \cite{Weil1,Weil2} that there exists
an automorphism $\phi: C \to C$, which induces (by functoriality)
either $\delta$ or $-\delta$ on $J(C)=X$. Replacing $\phi$ by
$\phi^4$ and taking into account that $\delta^3$ is the identity
automorphism of $X=J(C)$, we may and will assume that $\phi$ induces
$\delta$. Clearly, $\phi^3$ is the identity automorphism of $C$,
because it induces the identity map on $J(C)$ and $g>1$. The action
of $\phi$ on $C$ gives rise to the embedding
$$\mu_3 \hookrightarrow \Aut(C), \ \zeta_3 \mapsto \phi.$$

Let $P\in C$ be a fixed point of $\phi$. Then $\phi$ induces the automorphism of the
corresponding (one-dimensional) tangent space $T_P(C)$, which is multiplication by a complex
number $c_P$. Clearly, $c_P$ is a cubic root of unity.

\begin{lem}
\label{nondegenerate}
Every fixed point $P$ of $\phi$ is nondegenerate, i.e., $c_P\ne 1$.
\end{lem}

 \begin{proof}[Proof of Lemma \ref{nondegenerate}]    The result is  well-known. However, I failed to find a proper reference.

 Suppose that $c_P=1$.
 Let $\OC_P$ be the local ring at $P$ and $\m_P$ its maximal ideal. We write $\phi_{*}$ for
 the automorphism of $\OC_P$ induced by $\phi$. Clearly,    $\phi_{*}^3$    is the identity
 map. Since $\phi$ is {\sl not} the identity map, there
 are no $\phi_{*}$-invariant local parameters at $P$. Clearly,
  $\phi_{*}(\m_P)=\m_P,    \phi_{*}(\m_P^2)=\m_P^2$.
  Since      $T_P(C)$ is the dual of   $\m_P/\m_P^2$ and $c_p=1$, we conclude that
   $\phi_{*}$ induces the identity map on  $\m_P/\m_P^2$.
This implies that if $t\in \m_P$ is a local parameter at $t$ (i.e., its image $\bar{t}$ in
 $\m_P/\m_P^2$   is {\sl not} zero) then $t':=t+\phi_{*}(t)+  \phi_{*}^2(t)$ is
 $\phi_{*}$-invariant and its image in  $\m_P/\m_P^2$ equals $3\bar{t}\ne 0$. This implies
 that $t'$ ia a  $\phi_{*}$-invariant local parameter at $P$.      Contradiction.

 \end{proof}

 \begin{cor}
$D:=C/\mu_3$ is a smooth projective irreducible curve. The map $C\to D$ has degree $3$,
 its ramification points are exactly the images of  fixed points of $\phi$
 and all the ramification indices are $3$.
 \end{cor}

 \begin{lem}
 \label{line}
$D$ is biregularly isomorphic to the projective line.
 \end{lem}

 \begin{proof}[Proof of Lemma \ref{line}]
The map $C\to D$ induces, by Albanese functorialy, the surjective homomorphism
of the corresponding jacobians $J(C)\to J(D)$ that kills all the divisors classes
of the form $(Q)-(\phi(Q))$ ($Q\in C)$. This implies that it kills $(1-\delta)J(C)$.
On the other hand, $1-\delta:J(C)\to J(C)$ is, obviously, an isogeny. This implies that the
image of $J(C)$ in  $J(D)$ is zero and the surjectiveness implies that $J(D)=0$.
This means that the genus of $D$ is $0$.
 \end{proof}

 \begin{cor}
 \label{hur}
The number $h$ of fixed points of $\phi$ is $g+2$.
 \end{cor}

\begin{proof}[Proof of Corollary \ref{hur}]
Applying Hurwitz formula to $C\to D$, we get
$$2g-2=3\cdot (-2)+ 2\cdot h.$$
\end{proof}

\begin{lem}
\label{trace}
Let $\phi^{*}:\Omega^1(C) \to \Omega^1(C)$ be the automorphism of $\Omega^1(C)$ induced
by $\phi$ and $\tau$ its trace. Then
$$\tau=a \zeta_3+ b \zeta_3^{-1}.$$
\end{lem}

\begin{proof}[Proof of Lemma \ref{trace}]
Pick a $\phi$-invariant point $P_0$ and consider the regular map
$$\alpha:C \to J(C), Q \mapsto \cl((Q)-(P_0)).$$
It is well-known that $\alpha$ induces an isomorphism of complex vector spaces
$$\alpha^{*}:\Omega^1(X) \cong \Omega^1(C).$$
Clearly,
$$\phi^{*}=  \alpha^{*}\delta^{*}{\alpha^{*}}^{-1}$$
where $\delta^{*}:\Omega^1(J(C)) = \Omega^1(J(C))$
is the automorphism induced by $\delta$.
This implies that the traces of $\phi^{*}$ and $\delta^{*}$
do coincide. Now the very definition of $a$ and $b$ implies
that the trace of $\phi^{*}$ equals $a \zeta_3+ b \zeta_3^{-1}$.
\end{proof}

{\bf End of proof of Theorem \ref{nonjacob}}. Let $B$ be the set of fixed points of $\phi$.
 We know that $\#(B)=g+2$. By the holomorphic Lefschetz fixed point formula \cite[Th. 2]{AT},
\cite[Ch. 3, Sect. 4]{GH} (see also \cite[Sect. 12.2 and
12.5]{Milnor})
 applied to $\phi$,
 $$1-\bar{\tau}=\sum_{P\in B}\frac{1}{1-c_P}$$
 where $\bar{\tau}$ is the complex-conjugate of $\tau$. Recall that every $c_P$ is a (primitive)
  cubic root of unity and therefore
  $$\mid 1-c_P\mid =\sqrt{3},\ \mid\frac{1}{1-c_P}\mid =\frac{1}{\sqrt{3}}$$ and
  $$\mid 1-\bar{\tau}\mid \le \frac{g+2}{\sqrt{3}}.$$
  Now
  $$\mid 1-\bar{\tau}\mid^2=\frac{(a+b+2)^2+3(a-b)^2 }{4}=\frac{(g+2)^2+3(a-b)^2 }{4}.$$
  This implies that
  $$\frac{(g+2)^2}{{3}}\ge \frac{(g+2)^2+3(a-b)^2 }{4}.$$
  It follows that $(g+2)^2 \ge 9(a-b)^2$ and we are done.
       \end{proof}

\section{Cubic threefolds}
\label{threefolds}
Let $S:F(x_0,x_1,x_2,x_3)=0\subset \P^3$ be a smooth projective cubic surface over $\C$
 \cite{Manin}.
(In particular, $F$ is an irreducible homogeneous cubic polynomial in $x_0,x_1,x_2,x_3$
with complex coefficients.) Then the equation
$$y^3=F(x_0,x_1,x_2,x_3)$$
defines a smooth projective threefold $T\subset\P^4$ provided with the natural action
of $\mu_3$ that arises from multiplication of $y$ by cubic roots of unity  \cite{ACT}
(see also \cite{CT,MT}).     We have the $\mu_3$-invariant Hodge decomposition
$$\H^3(T,\C)=\H^3(T,\Z)\otimes\C=\H^{1,2}(T)\oplus \H^{2,1}(T)$$
and the $\mu_3$-invariant non-degenerate alternating intersection pairing
$$(,):   \H^3(T,\C) \times \H^3(T,\C) \to \C.$$
In addition,  both  $\H^{1,2}(T)$ and $\H^{2,1}(T)$ are $5$-dimensional isotropic subspaces
and $\mu_3$ acts on  $\H^{2,1}(T)$ with multiplicities   $(4,1)$, i.e.   $\zeta_3\in \mu_3$
acts as diagonalizable linear operator in  $\H^{2,1}(T)$ with eigenvalue $\zeta_3$ of multiplicity
 $4$  and   eigenvalue $\zeta_3^{-1}$ of multiplicity $1$
  (\cite[Sect. 5]{CT}, \cite[Sect. 2.2 and Lemma 2.6]{ACT}).  (The proof is based on
\cite[Th. 8.3 on p. 488]{G}; see also \cite[pp. 338--339]{GC}.)

Since   both  $\H^{1,2}(T)$ and $\H^{2,1}(T)$ are isotropic and the intersection pairing is
non-degenerate, its restriction    to    $\H^{1,2}(T)\times \H^{2,1}(T)$  gives rise to   the
non-degenerate  $\mu_3$-invariant  $\C$-bilinear pairing
$$(,): \H^{1,2}(T)\times \H^{2,1}(T) \to \C. \eqno(1)$$
It follows that   $\mu_3$ acts on  $\H^{1,2}(T)$  with
multiplicities   $(1,4)$. (This assertion also follows from the fact
that    $\H^{1,2}(T)$ is the complex-conjugate
 of  $\H^{2,1}(T)$.)
            In particular, the action of $\mu_3$ on  $\H^{1,2}(T)$  extends to the
 embedding
 $$\Z[\mu_3]\hookrightarrow    \End_{\C}(\H^{1,2}(T)).  \eqno(2)$$
\section{Intermediate jacobians}
\label{threef}
Let $(J(T),\theta_T)$ be the {\sl intermediate jacobian} of the cubic threefold $T$
 \cite[Sect. 3]{GC};
it is a principally polarized five-dimensional complex abelian variety. By functoriality,
$\mu_3$ acts on $J(T)$ and respects the principal polarization $\theta_T$.  As
complex torus,
$$J(T)=\H^{1,2}(T)/p(\H^3(T,\Z)),  \eqno(3)$$  where
$$p:  \H^3(T,\C)=\H^3(T,\Z)\otimes\C=\H^{1,2}(T)\oplus \H^{2,1}(T)\to \H^{1,2}(T)$$
is the projection map that kills   $\H^{2,1}(T)$.   The imaginary part of the Riemann
form of the polarization coincides with the intersection pairing on
$\H^3(T,\Z)\cong p(\H^3(T,\Z))$.

 It follows from (2) that the action
of $\mu_3$ on $J(T)$ extends to the embedding
$$\Z[\mu_3]\hookrightarrow    \End (J(T)).$$
Combining (1) and (3), we conclude that the $\mu_3$-modules
$\Omega^1(J(T))=\Hom_{\C}(\H^{1,2}(T),\C)$  and
$\H^{2,1}(T)$ are canonically isomorphic. Now the assertions of Sect. \ref{threefolds}
about multiplicities imply that  $\Z[\zeta_3]$ acts on
$\Omega^1(J(T))$ with multiplicities $(4,1)$.

Since $3 \times \mid 4-1\mid > 5+2$, it follows from Theorem \ref{nonjacob}
that $(J(T),\theta_T)$ is not isomorphic to the canonically polarized jacobian of a curve.
 Of course, this assertion was proven by completely different methods in \cite{GC}
for arbitrary smooth projective cubic threefolds.


\begin{thebibliography}{99}

\bibitem{ACT} D. Allcock, J.A. Carlson, D. Toledo, {\em The complex
hyperbolic geometry of the moduli space of cubic surfaces}.
J. Algebraic Geometry {\bf 11} (2002), 659--724.

\bibitem{AT}  M. F. Atiyah, R. Bott, {\em A Lefschetz fixed point formula
for elliptic differential operators}.
Bull. Amer. Math. Soc.  {\bf 72}  (1966), 245--250.

\bibitem{CT} J.A. Carlson, D. Toledo, {\em Discriminant complements
and kernels of monodromy representations}. Duke Math. J. {\bf 97} (1999), 621--648.

\bibitem{GC} C.H. Clemens, Ph. Griffiths, {\em The intermediate jacobian of
the cubic threefold}. Ann. of Math. (2) {\bf 95} (1972), 281--356.

%\bibitem{E} T. Ekedahl, {\em An effective version of Hilbert's irreducibility
%theorem}. In: S\'eminaire de Th\'eorie des Nombres, Paris 1988--89.
% Birkh\"auser, Boston, MA, 1990, pp. 241--249.

\bibitem{G} Ph. Griffiths, {\em On the periods of certain rational
 integrals}: I and II.  Ann. of Math. (2) {\bf 90} (1969), 460--541.


%\bibitem{Dolgachev} I. Dolgachev,  Topics in classical algebraic geometry, Part
%1, Chapter 10; available at
%http://www.math.lsa.umich.edu/~idolga/lecturenotes.html .



\bibitem{GH} Ph. Griffiths, J. Harris, Principles of algebraic geometry. John Wiley and Sons,
New York, 1978.

%\bibitem{Johnson} D.L. Johnson, {\em The group of formal power series under substitution}.
%J. Austral. Math. Soc. (Series A), {\bf 45} (1988), 296--302.

\bibitem{Manin} Yu. I. Manin, Cubic forms, second edition. North Holland, 1986.

\bibitem{MT} K. Matsumoto, T. Terasoma, {\em Theta constants associated
to cubic threefolds}. J. Algebraic Geometry {\bf 12} (2003), 741--775.

\bibitem{Milnor} J. Milnor, Dynamics in one complex variable. Vieweg, Braunschweig/Wiesbaden,
1999.

\bibitem{Weil1} A. Weil, {\em Zum Beweis des Torellischen Satzes}.
  G\"ott. Nachr. 1957, no. 2, pp. 33--53; \OE uvres, vol. III, [1957a].

\bibitem{Weil2} A. Weil, {\em Sur le th\'eor\`eme de Torelli}.
S\'eminaire Bourbaki {\bf 151} (Mai 1957).
  \end{thebibliography}
\end{document}